\documentclass[12pt]{article}
\usepackage{makeidx}
\usepackage{amsfonts}
\usepackage{amssymb}
\usepackage{amsmath}
\usepackage[numbers,sort&compress]{natbib}
\usepackage[colorlinks,citecolor=blue,urlcolor=blue]{hyperref}
\usepackage{graphicx,indentfirst,tabularx}
\usepackage{bm}
\usepackage{enumitem}
\usepackage{lmodern}
\usepackage{mathtools}
\usepackage{geometry}
\usepackage{booktabs}
\usepackage{array}
\usepackage{mathrsfs}

\newtheorem{theorem}{Theorem}

\newtheorem{corollary}[theorem]{Corollary}

\newtheorem{definition}[theorem]{Definition}

\newtheorem{lemma}[theorem]{Lemma}

\newtheorem{proposition}[theorem]{Proposition}

\newenvironment{proof}[1][Proof]{\noindent\textbf{#1.} }{\ \rule{0.5em}{0.5em}}
\input{tcilatex}
\begin{document}

\title{Zero detection via Cox random measures}
\author{Athanasios C. Micheas \thanks{%
Department of Statistics, University of Missouri, 146 Middlebush Hall,
Columbia, MO 65211-6100, USA, email: micheasa@missouri.edu}}
\maketitle

\begin{abstract}
Using point process theory we propose and investigate a probabilistic
identification of zeros of a deterministic real valued function. We
construct a random atomic driving measure for a Cox point process whose mass
at a candidate location depends explicitly on the value of the function at
that location. We study the properties of the constructed random measure,
including its exact asymptotic properties, its Laplace functional and
related random processes. We also give an important application involving
zero detection for Riemann's $\xi$ function.
\end{abstract}



\section{Introduction}

\label{sec:introduction}The problem of identifying the zeros of a real
valued function is fundamental in mathematical analysis, numerical
mathematics, probability, statistics and scientific computing. Given a
function $g:\mathcal{X}\rightarrow \mathbb{R}$, defined on a domain $%
\mathcal{X}\subseteq \mathbb{R}^{d}$, the objective is to identify the zero
set 
\begin{equation*}
\mathcal{Z}=\{x\in \mathcal{X}:g(x)=0\}.
\end{equation*}

Although classical numerical methods provide well established approaches for
locating individual roots, the identification of multiple zeros,
particularly in multidimensional domains, may present considerable
computational and analytical difficulties. These difficulties become more
pronounced when function evaluations are expensive, when the zero set has a
complicated geometrical structure, or when uncertainty is associated with
the available function evaluations.

Probabilistic approaches to zero finding have a substantial history. For
example, Strongin \cite{Strongin1972} developed a stochastic approach in
which the location of a root is described through a probability
distribution, and successive function evaluations are selected according to
the resulting probabilistic information. More generally, probabilistic
representations of unknown functions and their zeros have provided a basis
for statistical approaches to numerical analysis and optimization. Such
approaches motivate the investigation of alternative stochastic
constructions for identifying zero locations.

An important class of stochastic models is provided by point processes and
random measures. Point processes offer a natural mathematical framework for
representing random collections of locations, while their distributions can
be characterized through their associated counting measures and Laplace
functionals. In particular, Poisson point processes are characterized by an
intensity measure, which determines the expected number of points in
measurable subsets of the underlying domain. Cox processes, also known as
doubly stochastic Poisson processes, extend this construction by allowing
the intensity measure itself to be random. Conditional on a realization of
the driving random measure, a Cox process is an ordinary Poisson point
process. Standard treatments of these models and their mathematical
properties can be found in the texts by \cite{CoxIsham1980}, \cite%
{DaleyVereJones2003}, \cite{DaleyVereJones2008}, \cite{LastPenrose2017} and 
\cite{kallenberg2017random}.

For general point process theory on Poisson, Gibbs and Cox point processes
and their applications see the recent papers \cite{Micheas2014}, \cite%
{zhou2015spatio}, \cite{micheas2018sppmix}, \cite{Micheas2019}, \cite%
{ChenMicheasHolan2020}, \cite{TangLi2021}, \cite{chen2022hierarchical}, \cite%
{Baddeley2022}, \cite{WuMicheas2022}, \cite{chen2022hierarchical}, \cite%
{Baetal2023}, \cite{GoncalvesDias2023}, \cite{Micheas2025}, \cite%
{wang2026inference} and the references therein.

The use of Cox processes in spatial modeling has also been extended to
settings in which latent functions determine spatial intensity structures.
In particular, \cite{GoncalvesDias2023} developed a Bayesian methodology for
level-set Cox processes in which the partition of the spatial domain is
determined by a latent Gaussian process. This provides a relevant connection
between point process modeling and functionally defined spatial regions.

Motivated by these developments, this paper investigates a different use of
Cox point processes, namely, the probabilistic identification of zeros of a
deterministic real valued function. The central idea is to construct a
random atomic driving measure whose mass at a candidate location depends
explicitly on the value of the function at that location. In this way, the
stochastic intensity of the resulting Cox process incorporates information
about the zero structure of the function.

In this paper, we propose and study an intersection of three subjects:
numerical zero finding, stochastic point process theory, and asymptotic
classification using random measures. In particular, based on any real
function $g$, we propose a mathematical construction of a general Cox point
process which is able to identify the zeros set of $g$ and study its exact
and asymptotic properties.

The principal objective of this paper is to investigate the mathematical
consequences of this construction, with particular emphasis on the
asymptotic behavior of the driving random measures and the resulting Cox
point processes. The contributions can be summarized as follows.

First, we introduce a gamma driven atomic Cox point process whose intensity
measure incorporates the function dependent weight $\exp \{-K|g(x)|^{2}\}$.
The atomic component is essential for the identification of individual
candidate locations, since an intensity measure that is absolutely
continuous with respect to Lebesgue measure assigns zero mass to every
singleton.

Second, we derive the Laplace functional of the gamma driven random measure
and investigate its asymptotic behavior. We establish that the atomic mass
diverges in probability at candidate locations where the function vanishes,
whereas it converges to zero in probability at candidate locations where the
function is nonzero. This provides a measure theoretic description of the
separation induced by the proposed construction.

Third, we establish finite set statistical consistency of the detection
procedure. In particular, we derive an explicit upper bound for the
probability of at least one classification error and show that this
probability converges to zero as $K$ increases. Consequently, the complete
collection of candidate locations can be classified correctly with
probability approaching one.

Finally, we introduce a binary detection process, obtained by recording
whether the Cox process places at least one point at each candidate
location. We prove that this process converges in probability, as a finite
random measure, to the deterministic counting measure supported on the true
zeros within the candidate set. This result provides a process level
formulation of the asymptotic zero detection property.

An important distinction is that the proposed construction addresses the
probabilistic classification of a prescribed finite set of candidate
locations. It does not, by itself, provide an algorithm for generating
candidate locations or identifying an entire continuous zero set. These
questions require additional assumptions concerning the regularity of the
function, the geometry of its zero set, and the selection of candidate
locations, and these results will be forthcoming (see the discussion in the
concluding remarks).

The paper is organized as follows. Section 2 introduces the required point
process and random measure framework and defines the gamma driven Cox
construction, which is used in the construction of zero detection
probabilistic tests. Section 3 derives the Laplace functional and
establishes the asymptotic properties of the driving random measure.
Moreover, we investigate the detection probabilities and prove finite set
statistical consistency. Importantly, we establish convergence of the binary
detection process to the zero counting measure and discuss the selection of $%
K$ when it is kept fixed and finite. In Section 4, we apply the theory to
detection of zeros for Riemann's $\xi $ function, and therefore the
nontrivial zeros of Riemann's $\zeta $ function. Concluding remarks, with a
discussion of possible extensions, are given in the last section.

\section{Construction of zero detection tests via random measures}

The general definition for a point process $N$ and its distribution can be
given in terms of the induced probability measure. Indeed, let $N$ be a
measurable map from a probability space $(\Omega ,\mathcal{A},\mathbb{P})$
to the measurable space $(\mathcal{N}_{f}(\mathcal{X}),\mathcal{N})$, where $%
\mathcal{N}_{f}(\mathcal{X})$ is the collection of all locally finite
counting measures on $\mathcal{X}$, equipped with the $\sigma $-algebra
generated by the evaluation mappings. Formally, the evaluation $\sigma $%
-algebra $\mathcal{N}$ is the generated $\sigma $-algebra 
\begin{equation*}
\mathcal{N}=\sigma \left( \{v\in \mathcal{N}_{f}(\mathcal{X}%
):v(B)=n\}\right) ,
\end{equation*}%
for all $n=0,1,2,\dots ,$ and $B\in \mathcal{B}(\mathcal{X}),$ where $%
\mathcal{X}$ denotes the event space and $\mathcal{B}(\mathcal{X})$ the
Borel subsets of $\mathcal{X}.$ Then, the distribution of $N$ is the induced
probability measure given by%
\begin{equation*}
P_{N}(Y)=\mathbb{P}(N^{-1}(Y))=\mathbb{P}(N\in Y)=\mathbb{P}(\left\{ \omega
\in \Omega :N(\omega )\in Y\right\} ),
\end{equation*}%
for $Y\in \mathcal{N}.$ Note that by construction, the point process $N$ is
a random counting measure. This basic definition can be used to define any
point process $N$ and its distribution, and will be utilized in what
follows. In what follows, we use the notation $\overset{\mathbb{P}}{{%
\rightarrow }}$ to denote convergence in probability.

\subsection{Poisson point process}

Consider a region $\mathcal{W}\in \mathcal{B}(\mathcal{X})$, and suppose
that we observe $n$ points (events) $\varphi _{n}=\{\mathbf{s}%
_{i}\}_{i=1}^{n}=\mathbf{s}_{1:n},$ from a point process $N.$ In order to
model this collection of points we consider the inhomogeneous Poisson point
process (PPP) which requires two assumptions: first, the random variables $%
N(B)$, for measurable $B\subseteq \mathcal{W},$ are distributed as Poisson,
i.e., $N(B)\thicksim Pois(\Lambda (B)),$ and second, counts of events are
independent over any finite collection of disjoint regions.

The intensity measure $\Lambda (B)$ describes the average number of events
over $B$, and is defined by%
\begin{equation*}
\Lambda (B)=E\left[ N(B)\right] =\tint\limits_{B}\lambda (\mathbf{s})d\mu (%
\mathbf{s}),
\end{equation*}%
where $\lambda (\mathbf{s})$ is known as the intensity function. Note that $%
\Lambda $ need not itself be a counting measure, but it is a locally finite
measure. The intensity surface $\lambda (\mathbf{s})$ exists via an appeal
to the Radon-Nikodym theorem, provided that $\Lambda $ is an absolutely
continuous measure with respect to Lebesgue measure $\mu $ in $\mathcal{X}.$
The special case where $\lambda (\mathbf{s})=\lambda $ yields the
homogeneous Poisson process (HPP), with intensity $\lambda $ and mean
measure $\Lambda (B)=\lambda |B|\mathbf{,}$ where $|B|$ is the area of $B$.

The joint distribution of the independent events in $\varphi _{n},$ and the
number of events $N(\mathcal{W})=n$ over the window $\mathcal{W},$ is given
by%
\begin{equation}
f_{N}(\varphi _{n},n)=\frac{e^{-\Lambda (\mathcal{W})}}{n!}%
\tprod\limits_{i=1}^{n}\lambda (\mathbf{s}_{i})=\frac{1}{n!}%
e^{-\tint\limits_{\mathcal{W}}\lambda (\mathbf{s})d\mathbf{s}%
}\tprod\limits_{i=1}^{n}\lambda (\mathbf{s}_{i}),  \label{PoissonLik}
\end{equation}%
$n\geq 0.$ Note that the intensity function uniquely determines the PPP, and
the quantity $\lambda (\mathbf{s})d\mathbf{s}$ represents the approximate
expected number of points in an infinitesimal region of area $d\mathbf{s}$
around $\mathbf{s}$.

\subsection{Cox point process}

We give the general definition first.

\begin{definition}[Cox Point Process]
\label{DefCoxP}A point process $N_{C}$ is called a Cox (or doubly stochastic
Poisson) point processif, conditional on the realization of a random locally
finite measure $M,$ the process $N_{C}$ is an PPP with intensity measure $M. 
$ The distribution of $M$, denoted by $Q_{M}$, is called the driving
distribution. More precisely, the distribution $P_{N_{C}}$ of $N_{C}$ is
given by%
\begin{equation}
P_{N_{C}}(Y)=\tint\limits_{\mathcal{M}_{f}(\mathcal{X})}P_{\Lambda
}(Y)Q_{M}(d\Lambda ),  \label{EqDistrCoxPP1}
\end{equation}%
for all $Y\in \mathcal{N},$ where $P_{\Lambda }$ denotes the probability
distribution of a Poisson point process with deterministic intensity measure 
$\Lambda $, and $\mathcal{M}_{f}(\mathcal{X})$ denotes the space of all
locally finite measures.
\end{definition}

The Cox point process (CPP) can be viewed as a hierarchical model where at
the second stage we have a realization of an intensity measure $M=\Lambda $
from some distribution $Q_{M}$, and at the first (data) stage we observe the
data from an PPP $N$ with intensity measure $\Lambda .$ Therefore, $%
M\thicksim Q_{M}$ and $N_{C}|M=\Lambda \thicksim P_{\Lambda },$ define a
hierarchical mixture distribution, with the unconditional (or marginal)
distribution being the distribution of the Cox process $N_{C}.$ The
distribution $Q_{M}$ is known as the driving distribution.

Now assume that the driving measure is a random field%
\begin{equation}
M=\{\lambda _{\mathbf{\theta }}(\mathbf{s})\geq 0:\mathbf{\theta }\thicksim
\pi (\mathbf{\theta }|\mathbf{\psi }),\mathbf{s}\in \mathcal{W},\mathbf{%
\theta }\in \mathbf{\Theta }\},  \label{RandomField}
\end{equation}%
with%
\begin{equation*}
M(B)=\tint\limits_{B}\lambda _{\mathbf{\theta }}(\mathbf{s})d\mu (\mathbf{s}%
),
\end{equation*}%
for all $B\in \mathcal{B}(\mathcal{W}),$ and $\lambda _{\mathbf{\theta }}(%
\mathbf{s})$ must be nonnegative and integrable over bounded regions in
order to define a valid locally finite intensity measure.

The source of randomness in $M$ are the parameters $\mathbf{\theta }$ which
are modeled according to some density providing the source of randomness
(hereafter called the source or driving density) $\pi (\mathbf{\theta }|%
\mathbf{\psi }),$ for some hyperparameters $\mathbf{\psi }$ that require
estimation$.$

In order to sample from this Cox point process, we first sample $\mathbf{%
\theta }_{0}\thicksim \pi (\mathbf{\theta }|\mathbf{\psi }),$ so that the
realized intensity function is $\lambda _{\mathbf{\theta }_{0}}(\mathbf{s}),$
$\mathbf{s}\in \mathcal{X}$ (a realization of the random field). Then
sampling an PPP with intensity function $\lambda _{\mathbf{\theta }_{0}}(%
\mathbf{s}),$ gives a realization from a Cox point process $N_{C}$ with
driving measure $M$ given by (\ref{RandomField}).

Given the realization of the intensity function $\lambda _{\mathbf{\theta }}(%
\mathbf{s}),$ the density of the PPP $N$ is given by (\ref{PoissonLik}). In
view of (\ref{EqDistrCoxPP1}), we can write the density $f_{N_{C}}(\varphi
_{n},n)$ of $N_{C}$ as the mixture distribution%
\begin{equation}
f_{N_{C}}(\varphi _{n},n|\mathbf{\psi })=\tint\limits_{\Theta }f_{N}(\varphi
_{n},n|\mathbf{\theta })\pi (\mathbf{\theta }|\mathbf{\psi })d\mathbf{\theta 
},  \label{EqDistrCoxPP}
\end{equation}%
for all point patterns $\varphi _{n}$ with $n=0,1,\ldots ,$ and as a result
the likelihood associated with one observed realization (point pattern) $%
\varphi _{n}$ of $N_{C}$ is given by%
\begin{equation}
f_{N_{C}}(\varphi _{n},n|\mathbf{\psi })=\frac{1}{n!}\tint\limits_{\mathbf{%
\Theta }}\exp \left\{ -\tint\limits_{\mathcal{W}}\lambda _{\mathbf{\theta }}(%
\mathbf{s})d\mathbf{s}\right\} \left[ \tprod\limits_{i=1}^{n}\lambda _{%
\mathbf{\theta }}(\mathbf{s}_{i})\right] \pi (\mathbf{\theta }|\mathbf{\psi }%
)d\mathbf{\theta }.  \label{CPLikelihood}
\end{equation}

\subsection{Exact tests for a zero location}

We are interested in identifying the zeros of a real function $g:\mathcal{X}%
\rightarrow \mathbb{R},$\ based on an appropriate construction of a CPP that
depends on $g$. Define the deterministic zero set 
\begin{equation*}
\mathcal{Z}:=\{x\in \mathcal{X}:g(x)=0\},
\end{equation*}%
and consider the CPP $N_{\mathcal{Z}}$ whose realizations consist of points
belonging to $\mathcal{Z}.$ The Cox point process will act as a
probabilistic zero detection mechanism: given a test location $x_{0}$, the
process should determine whether $g(x_{0})=0$, without being explicitly told
whether $x_{0}$ belongs to the zero set. Ideally, we want the following
equivalence%
\begin{equation}
\mathbf{P}(N_{\mathcal{Z}}(\left\{ x_{0}\right\} )\geq 1)>0\Leftrightarrow
g(x_{0})=0,  \label{cond1}
\end{equation}%
or equivalently%
\begin{equation}
\mathbf{P}(N_{\mathcal{Z}}(\left\{ x_{0}\right\} )\geq 1)=0\Leftrightarrow
g(x_{0})\neq 0,  \label{cond2}
\end{equation}%
This would give us the desired exact zero detection; the process identifies
a point $x_{0}$ if and only if $g(x_{0})=0$.

However, an ordinary Cox process with an intensity measure absolutely
continuous with respect to Lebesgue measure cannot satisfy this condition,
since $M(\{x_{0}\})=0,$ regardless of whether $g(x_{0})=0$. Therefore, the
dominating measure $\mu $ is required to have non-atomic (e.g., Lebesgue)
and atomic parts (e.g., Dirac). In particular, to identify individual zeros,
the driving random measure must be capable of assigning positive mass to
individual locations.

For $K>0$, define the intensity function of a Poisson point process $N_{K,%
\mathbf{\theta }}$ on $\mathcal{X}$ by 
\begin{equation}
\lambda _{K,\mathbf{\theta }}(x)=e^{-K|g(x)|^{2}},  \label{intens1}
\end{equation}%
for $x\in \mathcal{X},$ where $\mathbf{\theta }$ is a realization of the
driving density $\pi (\mathbf{\theta }|\mathbf{\psi }).$ We will construct
the driving measure of $N_{\mathcal{Z}}$ in such a way as to have $N_{%
\mathcal{Z}}$ satisfy equation (\ref{cond1}) or equivalently, equation (\ref%
{cond2}).

We begin with an exact result.

\begin{lemma}[Exact zero detection by a Cox point process]
\label{lem:exact_zero_detection} Let $\mathcal{X}$ be a measurable space and
let $g:\mathcal{X}\rightarrow \mathbb{R}$ be a measurable function with zero
set 
\begin{equation*}
\mathcal{Z}:=\{x\in \mathcal{X}:g(x)=0\}.
\end{equation*}%
Let $D=\{x_{j}:j\in J\}\subseteq \mathcal{X}$ be a countable set of
candidate locations, where $J$ is finite or countable. Consider a Cox point
process $N_{\mathcal{Z}}$ with random driving measure 
\begin{equation*}
M_{\theta }=\sum_{j\in J}A_{j,\theta }1(x_{j}\in \mathcal{Z})\delta
_{x_{j}}+M_{\theta }^{(c)},
\end{equation*}%
where $\delta _{x_{j}}$ denotes the Dirac measure at $x_{j}$, $A_{j,\theta }$
are nonnegative random variables, each interpreted as random mass assigned
to location $x_{j},$ satisfying 
\begin{equation*}
\mathbb{P}(A_{j,\theta }>0)>0,
\end{equation*}%
$j\in J,$ and $M_{\theta }^{(c)}$ is a diffuse random measure, that is, 
\begin{equation*}
M_{\theta }^{(c)}(\{x\})=0\quad \text{almost surely for every }x\in \mathcal{%
X}.
\end{equation*}%
Assume that $M_{\theta }$ is locally finite almost surely.

Conditional on $M_{\theta }$, let $N_{\mathcal{Z}}$ be a Poisson point
process with intensity measure $M_{\theta }$. Then, for every $x_{j}\in D$, 
\begin{equation*}
{\mathbb{P}}\left( N_{\mathcal{Z}}{(\{x_{j}\})\geq 1}\right)
>0\Leftrightarrow g(x_{j})=0{,}
\end{equation*}%
or equivalently, we have%
\begin{equation*}
{\mathbb{P}}\left( N_{\mathcal{Z}}{(\{x_{j}\})\geq 1}\right)
=0\Leftrightarrow g(x_{j})\neq 0{.}
\end{equation*}%
Consequently, the Cox point process provides exact probabilistic zero
detection on the candidate set $D$.
\end{lemma}

\begin{proof}
By the defining property of a Cox point process, conditional on the
realization of the driving measure $M_{\theta }$, 
\begin{equation*}
N_{\mathcal{Z}}(\{x_{j}\})\mid M_{\theta }\sim \func{Pois}\left( M_{\theta
}(\{x_{j}\})\right) .
\end{equation*}%
Since the diffuse component assigns zero mass to singletons, we have 
\begin{equation*}
M_{\theta }(\{x_{j}\})=A_{j,\theta }1(x_{j}\in \mathcal{Z}).
\end{equation*}%
Therefore, integrating out the random driving measure, 
\begin{equation*}
\mathbb{P}\left( N_{\mathcal{Z}}(\{x_{j}\})\geq 1\right) =1-\mathbb{E}\left[
\exp \{-M_{\theta }(\{x_{j}\})\}\right] =1-\mathbb{E}\left[ \exp \left\{
-A_{j,\theta }1(x_{j}\in \mathcal{Z})\right\} \right] .
\end{equation*}%
If $g(x_{j})\neq 0$, then $x_{j}\notin \mathcal{Z}$, and hence 
\begin{equation*}
\mathbb{P}\left( N_{\mathcal{Z}}(\{x_{j}\})\geq 1\right) =0.
\end{equation*}%
On the other hand, if $g(x_{j})=0$, then 
\begin{equation*}
\mathbb{P}\left( N_{\mathcal{Z}}(\{x_{j}\})\geq 1\right) =1-\mathbb{E}\left[
e^{-A_{j,\theta }}\right] .
\end{equation*}%
Since $A_{j,\theta }\geq 0$ and $\mathbb{P}(A_{j,\theta }>0)>0$, it follows
that 
\begin{equation*}
1-\mathbb{E}\left[ e^{-A_{j,\theta }}\right] >0.
\end{equation*}%
This establishes both equivalences.
\end{proof}

Next we take the detection process to the limit with respect to $K$.

\begin{lemma}[Asymptotically exact zero detection by a Cox point process]
\label{lem:asymptotic_zero_detection}

Let $\mathcal{X}$ be a measurable space and let $g:\mathcal{X}\rightarrow 
\mathbb{R}$ be a measurable function. Consider a countable set of candidate
locations 
\begin{equation*}
D=\{x_{j}:j\in J\}\subseteq \mathcal{X},
\end{equation*}%
where $J$ is finite or countable.

For $K>0$, define the function 
\begin{equation*}
\lambda _{K}(x)=e^{-K|g(x)|^{2}},\qquad x\in \mathcal{X}.
\end{equation*}%
Let $\{C_{j,\theta }:j\in J\}$ be a collection of strictly positive random
variables satisfying 
\begin{equation*}
C_{j,\theta }>0\quad \text{almost surely},\qquad \mathbb{E}[C_{j,\theta
}]<\infty .
\end{equation*}%
Define the random atomic measure 
\begin{equation*}
M_{K,\theta }^{(a)}=\sum_{j\in J}C_{j,\theta }e^{-K|g(x_{j})|^{2}}\delta
_{x_{j}}.
\end{equation*}%
Furthermore, let $W_{\theta }$ be a diffuse random measure on $\mathcal{X}$,
with locally finite realizations, and define the diffuse component 
\begin{equation*}
M_{K,\theta }^{(c)}(B)=\tint\limits_{B}e^{-K|g(x)|^{2}}W_{\theta }(dx),
\end{equation*}%
$B\in \mathcal{B}(X).$ Assume that the resulting random measure 
\begin{equation*}
M_{K,\theta }=M_{K,\theta }^{(a)}+M_{K,\theta }^{(c)}
\end{equation*}%
is locally finite almost surely.

Conditional on $M_{K,\theta }$, let $N_{K}$ be a Poisson point process with
intensity measure $M_{K,\theta }$. Then, for every candidate location $%
x_{j}\in D$, we have 
\begin{equation*}
{\lim_{K\rightarrow +\infty }\mathbb{P}}\left( {N_{K}(\{x_{j}\})\geq 1}%
\right) {=%
\begin{cases}
1-\mathbb{E}[e^{-C_{j,\theta }}]>0, & g(x_{j})=0, \\[2mm] 
0, & g(x_{j})\neq 0,%
\end{cases}%
}
\end{equation*}%
and consequently 
\begin{equation*}
{\lim_{K\rightarrow +\infty }\mathbb{P}}\left( {N_{K}(\{x_{j}\})\geq 1}%
\right) {>0\Leftrightarrow g(x_{j})=0.}
\end{equation*}%
Thus, the proposed Cox point process provides asymptotically exact
probabilistic zero detection on the candidate set $D$, without explicitly
using the zero-set indicator in the construction of its driving measure.
\end{lemma}

\begin{proof}
By the definition of a Cox point process, conditional on the driving measure 
$M_{K,\theta }$, we have 
\begin{equation*}
N_{K}(\{x_{j}\})\mid M_{K,\theta }\sim \func{Pois}\left( M_{K,\theta
}(\{x_{j}\})\right) .
\end{equation*}%
Since $W_{\theta }$ is diffuse, its contribution to the mass of a singleton
is zero. Therefore, 
\begin{equation*}
M_{K,\theta }(\{x_{j}\})=C_{j,\theta }e^{-K|g(x_{j})|^{2}}.
\end{equation*}%
Integrating out the random driving measure gives 
\begin{equation*}
\mathbb{P}\left( N_{K}(\{x_{j}\})\geq 1\right) =1-\mathbb{E}\left[ \exp
\left\{ -C_{j,\theta }e^{-K|g(x_{j})|^{2}}\right\} \right] .
\end{equation*}%
We distinguish two cases.

\medskip \noindent \textbf{Case 1: $g(x_j)=0$.}

In this case, 
\begin{equation*}
e^{-K|g(x_{j})|^{2}}=1
\end{equation*}%
for every $K>0$. Consequently, 
\begin{equation*}
\mathbb{P}\left( N_{K}(\{x_{j}\})\geq 1\right) =1-\mathbb{E}[e^{-C_{j,\theta
}}]>0,
\end{equation*}%
where strict positivity follows from $C_{j,\theta }>0$ almost surely.

\medskip \noindent \textbf{Case 2: $g(x_j)\neq0$.}

Since $|g(x_{j})|^{2}>0$, we have 
\begin{equation*}
\lim_{K\rightarrow +\infty }e^{-K|g(x_{j})|^{2}}=0.
\end{equation*}%
It follows that 
\begin{equation*}
\lim_{K\rightarrow +\infty }\exp \left\{ -C_{j,\theta
}e^{-K|g(x_{j})|^{2}}\right\} =1\quad \text{almost surely}.
\end{equation*}%
Moreover, the integrand is bounded between zero and one. Hence, by the
dominated convergence theorem, 
\begin{equation*}
\lim_{K\rightarrow +\infty }\mathbb{P}\left( N_{K}(\{x_{j}\})\geq 1\right)
=1-\mathbb{E}[1]=0.
\end{equation*}%
Combining the two cases establishes the result.
\end{proof}

Its purpose is to describe the spatial distribution of points near the zero
set. The atomic component provides the pointwise detection mechanism.

One technical qualification; the latter lemma establishes convergence of the
detection probabilities at each fixed candidate location. It does not
establish uniform convergence over an arbitrary infinite candidate set.
Also, for finite $K$, false detections remain possible. The construction is
therefore an asymptotically exact detection mechanism, not an exact finite $%
K $ test.

\subsection{Zero detection by a gamma driven Cox point process}

Next we consider a specific model for the underlying source of randomness of
the random field of the Cox process.

\begin{lemma}[gamma driven atomic masses and zero detection probability]
\label{lem:gamma_zero_detection}

Let $g:\mathcal{X}\rightarrow \mathbb{R}$ be a measurable function and let $%
D=\{x_{j}:j\in J\}$ be a countable set of candidate locations.

For each $x_{j}\in D$, let the random atomic mass $C_{j}$ independently
follow a Gamma distribution with shape parameter $a_{j}>0$ and rate
parameter $b_{j}>0$, that is, 
\begin{equation*}
C_{j}\sim \func{Gamma}(a_{j},b_{j}),
\end{equation*}%
with density 
\begin{equation*}
\pi (c_{j}\mid a_{j},b_{j})=\frac{b_{j}^{a_{j}}}{\Gamma (a_{j})}%
c_{j}^{a_{j}-1}e^{-b_{j}c_{j}},\qquad c_{j}>0.
\end{equation*}%
For $K>0$, define the random atomic measure 
\begin{equation*}
M_{K}^{(a)}=\sum_{j\in J}C_{j}e^{-K|g(x_{j})|^{2}}\delta _{x_{j}}.
\end{equation*}%
Let $M_{K}^{(c)}$ be a diffuse random measure and define 
\begin{equation*}
M_{K}=M_{K}^{(a)}+M_{K}^{(c)}.
\end{equation*}%
Assume that $M_{K}$ is locally finite almost surely. Conditional on $M_{K}$,
let $N_{K}$ be a Poisson point process with intensity measure $M_{K}$.

Then, for every candidate location $x_{j}\in D$, the probability of
detecting at least one point at $x_{j}$ is given by 
\begin{equation*}
{p_{K}(x_{j})=1-\left( \frac{b_{j}}{b_{j}+e^{-K|g(x_{j})|^{2}}}\right)
^{a_{j}}.}
\end{equation*}%
Moreover%
\begin{equation*}
\lim_{K\rightarrow +\infty }p_{K}(x_{j})=%
\begin{cases}
1-\left( \dfrac{b_{j}}{b_{j}+1}\right) ^{a_{j}}, & g(x_{j})=0, \\[3mm] 
0, & g(x_{j})\neq 0,%
\end{cases}%
\end{equation*}%
and therefore 
\begin{equation*}
\lim_{K\rightarrow +\infty }p_{K}(x_{j})>0\Leftrightarrow g(x_{j})=0.
\end{equation*}
\end{lemma}

\begin{proof}
Since $M_{K}^{(c)}$ is diffuse, its mass at any singleton is zero.
Therefore, 
\begin{equation*}
M_{K}(\{x_{j}\})=C_{j}e^{-K|g(x_{j})|^{2}}.
\end{equation*}%
By the defining property of a Cox point process, 
\begin{equation*}
N_{K}(\{x_{j}\})\mid M_{K}\sim \func{Pois}(M_{K}(\{x_{j}\})).
\end{equation*}%
Taking expectations with respect to the random driving measure, we obtain 
\begin{equation*}
p_{K}(x_{j})=\mathbb{P}(N_{K}(\{x_{j}\})\geq 1)=1-\mathbb{E}\left[
e^{-C_{j}e^{-K|g(x_{j})|^{2}}}\right] .
\end{equation*}%
Using the Laplace transform of a Gamma random variable, 
\begin{equation*}
\mathbb{E}[e^{-tC_{j}}]=\left( \frac{b_{j}}{b_{j}+t}\right) ^{a_{j}},\qquad
t\geq 0,
\end{equation*}%
and substituting 
\begin{equation*}
t=e^{-K|g(x_{j})|^{2}},
\end{equation*}%
we obtain 
\begin{equation*}
p_{K}(x_{j})=1-\left( \frac{b_{j}}{b_{j}+e^{-K|g(x_{j})|^{2}}}\right)
^{a_{j}}.
\end{equation*}%
If $g(x_{j})=0$, then $e^{-K|g(x_{j})|^{2}}=1$, and hence 
\begin{equation*}
p_{K}(x_{j})=1-\left( \frac{b_{j}}{b_{j}+1}\right) ^{a_{j}}>0.
\end{equation*}%
If $g(x_{j})\neq 0$, then 
\begin{equation*}
\lim_{K\rightarrow \infty }e^{-K|g(x_{j})|^{2}}=0,
\end{equation*}%
which implies 
\begin{equation*}
\lim_{K\rightarrow \infty }p_{K}(x_{j})=0.
\end{equation*}%
This establishes the result.
\end{proof}

The latter lemma reveals an important underlying structure of the detection
test; the limiting probability of detecting a true zero is controlled by $%
a_{j}${}and $b_{j}$. In particular, increasing $a_{j}$, while keeping $b_{j}$
fixed, increases detection probability, whilst increasing $b_{j}$, while
keeping $a_{j}$ fixed, decreases detection probability. Moreover, the ratio $%
a_{j}/b_{j}$ determines the expected atomic mass.

The next result follows trivially.

\begin{corollary}[Marginal distribution of the atomic counts]
\label{cor:negative_binomial}

Under the assumptions of Lemma~\ref{lem:gamma_zero_detection}, the counts at
distinct candidate locations are independent after integrating out the
independent Gamma atomic masses.

In particular, we have 
\begin{equation*}
N_{K}(\{x_{j}\})\sim \func{NegBin}\left( a_{j},\frac{b_{j}}{b_{j}+q_{j}(K)}%
\right) ,
\end{equation*}%
where the second parameter denotes the success probability.

The corresponding probability mass function is 
\begin{equation*}
\mathbb{P}(N_{K}(\{x_{j}\})=n_{j})=\frac{\Gamma (a_{j}+n_{j})}{\Gamma
(a_{j})n_{j}!}\left( \frac{b_{j}}{b_{j}+q_{j}(K)}\right) ^{a_{j}}\left( 
\frac{q_{j}(K)}{b_{j}+q_{j}(K)}\right) ^{n_{j}},
\end{equation*}%
for $n_{j}=0,1,2,\ldots $.
\end{corollary}

This is a particularly useful result from a statistical modeling point of
view, because it shows us that we can estimate $a_{j}$, $b_{j}$, and
potentially $K$ from observed point counts, using for example maximum
likelihood or Bayesian inference. Moreover, even as $K\rightarrow +\infty $,
a true zero can be missed in an individual realization. This is not a
mathematical inconsistency; it simply means that the model distinguishes
between detection capability, where a true zero has positive detection
probability and certain detection, where a true zero is detected with
probability one.

When the second property is desired, one possible extension is to let the
Gamma shape parameter increase with $K$, while controlling the rate
parameter appropriately. This requires a separate analysis to ensure that
the false positive probability still converges to zero. We collect this
result in the following.

\begin{lemma}[Asymptotically certain zero detection]
\label{lem:certain_zero_detection}

Let $g:\mathcal{X}\rightarrow \mathbb{R}$ be a measurable function and let $%
D=\{x_{j}:j\in J\}$ be a countable set of candidate locations.

For each $K>0$, let the random atomic masses be independently distributed
according to 
\begin{equation*}
C_{j,K}\sim\func{Gamma}(K,1), \qquad j\in J,
\end{equation*}
where the second parameter denotes the rate.

Define the random atomic measure 
\begin{equation*}
M_{K}^{(a)}=\sum_{j\in J}C_{j,K}e^{-K|g(x_{j})|^{2}}\delta _{x_{j}}.
\end{equation*}%
Let $M_{K}^{(c)}$ be a diffuse random measure such that 
\begin{equation*}
M_{K}=M_{K}^{(a)}+M_{K}^{(c)}
\end{equation*}%
is locally finite almost surely.

Conditional on $M_{K}$, let $N_{K}$ be a Poisson point process with
intensity measure $M_{K}$. Then, for every candidate location $x_{j}\in D$, 
\begin{equation*}
\mathbb{P}\left( {N_{K}(\{x_{j}\})\geq 1}\right) =1-\left(
1+e^{-K|g(x_{j})|^{2}}\right) ^{-K}.
\end{equation*}%
Moreover, we have 
\begin{equation*}
\lim_{K\rightarrow +\infty }\mathbb{P}\left( {N_{K}(\{x_{j}\})\geq 1}\right)
=1(g(x_{j})=0),
\end{equation*}%
or equivalently, we have%
\begin{align*}
g(x_{j})& =0\Rightarrow \mathbb{P}\left( N_{K}(\{x_{j}\})\geq 1\right)
\rightarrow 1, \\
g(x_{j})& \neq 0\Rightarrow \mathbb{P}\left( N_{K}(\{x_{j}\})\geq 1\right)
\rightarrow 0.
\end{align*}%
Thus, the proposed Cox point process provides asymptotically certain zero
detection on the candidate set $D$.
\end{lemma}

\begin{proof}
Since $M_{K}^{(c)}$ is diffuse, its mass at a singleton is zero. Therefore, 
\begin{equation*}
M_{K}(\{x_{j}\})=C_{j,K}e^{-K|g(x_{j})|^{2}}.
\end{equation*}%
By the defining property of a Cox point process, 
\begin{equation*}
N_{K}(\{x_{j}\})\mid M_{K}\sim \func{Pois}(M_{K}(\{x_{j}\})).
\end{equation*}%
Consequently, 
\begin{equation*}
\mathbb{P}\left( N_{K}(\{x_{j}\})\geq 1\right) =1-\mathbb{E}\left[
e^{-C_{j,K}e^{-K|g(x_{j})|^{2}}}\right] .
\end{equation*}%
Using the Laplace transform of the $\func{Gamma}(K,1)$ distribution, 
\begin{equation*}
\mathbb{E}[e^{-tC_{j,K}}]=(1+t)^{-K},\qquad t\geq 0,
\end{equation*}%
we obtain 
\begin{equation*}
\mathbb{P}\left( N_{K}(\{x_{j}\})\geq 1\right) =1-\left(
1+e^{-K|g(x_{j})|^{2}}\right) ^{-K}.
\end{equation*}%
We now consider two cases.

\medskip \noindent \textbf{Case 1: $g(x_j)=0$.}

In this case, 
\begin{equation*}
\mathbb{P}\left( N_{K}(\{x_{j}\})\geq 1\right) =1-(1+1)^{-K}=1-2^{-K},
\end{equation*}%
and therefore 
\begin{equation*}
\lim_{K\rightarrow +\infty }\mathbb{P}\left( N_{K}(\{x_{j}\})\geq 1\right)
=1.
\end{equation*}

\medskip \noindent \textbf{Case 2: $g(x_j)\neq0$.}

Let $c=|g(x_{j})|^{2}>0$. Then 
\begin{equation*}
0\leq \mathbb{P}\left( N_{K}(\{x_{j}\})\geq 1\right) =1-\left(
1+e^{-Kc}\right) ^{-K}\leq Ke^{-Kc},
\end{equation*}%
where we have used the inequality 
\begin{equation*}
1-(1+t)^{-K}\leq Kt,\qquad t\geq 0.
\end{equation*}%
Since 
\begin{equation*}
\lim_{K\rightarrow \infty }Ke^{-Kc}=0,
\end{equation*}%
the probability is bounded above and below by zero and therefore 
\begin{equation*}
\lim_{K\rightarrow +\infty }\mathbb{P}\left( N_{K}(\{x_{j}\})\geq 1\right)
=0.
\end{equation*}%
Combining the two cases establishes the result.
\end{proof}

\subsection{Statistical interpretation and detection properties}

\label{subsec:detection_properties}The proposed gamma driven Cox point
process provides an interesting interpretation of the parameter $K$ as a
zero detection concentration parameter. In particular, increasing $K$
simultaneously increases the probability of detecting a true zero and
decreases the probability of falsely identifying a nonzero location as a
zero.

Indeed, for a candidate location $x_{j}\in D$, the probability of failing to
detect a true zero is given by 
\begin{equation*}
\mathbb{P}\left( N_{K}(\{x_{j}\})=0\right) =\left(
1+e^{-K|g(x_{j})|^{2}}\right) ^{-K}.
\end{equation*}%
In particular, if $g(x_{j})=0$, this expression reduces to 
\begin{equation*}
\mathbb{P}\left( N_{K}(\{x_{j}\})=0\right) =2^{-K},
\end{equation*}%
which converges to zero as $K\rightarrow +\infty $. Consequently, the
probability of detecting a true zero converges to one.

On the other hand, if $g(x_{j})\neq 0$, the probability of falsely
identifying $x_{j}$ as a zero satisfies 
\begin{equation*}
\mathbb{P}\left( N_{K}(\{x_{j}\})\geq 1\right) \leq Ke^{-K|g(x_{j})|^{2}},
\end{equation*}%
which converges to zero as $K\rightarrow +\infty $.

Thus, increasing $K$ provides a simultaneous reduction in the probability of
false negatives and false positives at every fixed candidate location. These
properties establish $K$ as a parameter controlling the trade off between
the concentration of the driving measure around the zero set and the
accuracy of the resulting zero detection mechanism.

\subsection{Further considerations and limitations}

\label{subsec:limitations}

Although the proposed construction provides asymptotically certain zero
detection at every fixed candidate location, some additional considerations
are necessary.

First, the expected atomic mass at a candidate location is given by 
\begin{equation*}
\mathbb{E}[C_{j,K}]=K,
\end{equation*}%
whereas its variance is 
\begin{equation*}
\func{Var}(C_{j,K})=K.
\end{equation*}%
Consequently, the expected atomic mass at a true zero increases without
bound as $K\rightarrow +\infty $. This increase is responsible for the
convergence of the detection probability to one.

However, this also implies that the sequence of driving measures need not
converge to a locally finite random measure. The construction should
therefore be interpreted as a sequence of locally finite Cox processes
indexed by $K$, rather than necessarily as a single limiting Cox process.

Second, the detection result established in Lemma~\ref%
{lem:certain_zero_detection} is pointwise in the candidate location. That
is, for every fixed $x_{j}\in D$, the detection probability converges to $%
1(g(x_{j})=0)$. This does not automatically establish uniform convergence
over an infinite candidate set or a continuous domain.

For a finite candidate set $D=\{x_{1},\ldots ,x_{m}\}$, however, the
pointwise convergence also implies simultaneous convergence in probability.
In particular, 
\begin{equation*}
\mathbb{P}\left( 
\begin{array}{c}
\text{at least one detection error} \\ 
\text{occurs in }D%
\end{array}%
\right) \leq \sum_{j:g(x_{j})=0}2^{-K}+\sum_{j:g(x_{j})\neq
0}Ke^{-K|g(x_{j})|^{2}},
\end{equation*}%
which converges to zero as $K\rightarrow +\infty $.

Therefore, on any finite candidate set, the proposed Cox point process
identifies the zero and nonzero locations simultaneously with probability
approaching one.

Finally, it is important to distinguish the finite $K$ model from its
asymptotic behavior. For any finite $K>0$, the detection probability at a
nonzero location remains strictly positive, and hence false positive
detections cannot be completely excluded. Exact zero detection is obtained
only in the asymptotic sense as $K$ increases without bound.

These observations motivate the use of the proposed gamma driven Cox process
as a probabilistic zero detection model, in which the parameter $K$ controls
the accuracy of the detection procedure.

\section{Laplace functional and the binary detection process}

First, we derive the Laplace functional along with asymptotic behavior of
the random atomic measure in the following.

\begin{lemma}[Asymptotic behavior of gamma driven atomic random measures]
\label{lem:random_measure_asymptotics}

Let $\mathcal{X}$ be a locally compact Polish space and let 
\begin{equation*}
D=\{x_{1},\ldots ,x_{m}\}\subset \mathcal{X}
\end{equation*}%
be a finite set of distinct candidate locations. Let $g:\mathcal{X}%
\rightarrow \mathbb{R}$ be measurable. For each $K>0$, let 
\begin{equation*}
C_{j,K}\sim \func{Gamma}(K,1),\qquad j=1,\ldots ,m,
\end{equation*}%
be independent random variables, where the second parameter denotes the
rate. Define the random atomic measure 
\begin{equation*}
M_{K}^{(a)}=\sum_{j=1}^{m}C_{j,K}e^{-K|g(x_{j})|^{2}}\delta _{x_{j}}.
\end{equation*}%
Let 
\begin{equation*}
\mathcal{Z}_{D}=\{x_{j}\in D:g(x_{j})=0\}
\end{equation*}%
denote the zeros of $g$ in $D$, and let 
\begin{equation*}
D_{+}=D\setminus \mathcal{Z}_{D}
\end{equation*}%
denote the nonzero candidate locations.

Then the following statements hold.

\medskip \noindent \textbf{(i) Laplace functional.}

For every nonnegative continuous function $f\in C_{c}(\mathcal{X})$, the
Laplace functional of $M_{K}^{(a)}$ is given by 
\begin{equation*}
{L_{K}(f)=\mathbb{E}\left[ e^{-\int\limits_{\mathcal{X}}f(x)M_{K}^{(a)}(dx)}%
\right] =\prod_{j=1}^{m}\left( 1+f(x_{j})e^{-K|g(x_{j})|^{2}}\right) ^{-K}.}
\end{equation*}

\medskip \noindent \textbf{(ii) Behavior at true zeros.}

For every $x_{j}\in \mathcal{Z}_{D}$, 
\begin{equation*}
C_{j,K}\overset{\mathbb{P}}{{\rightarrow }}+\infty .
\end{equation*}%
Consequently, 
\begin{equation*}
M_{K}^{(a)}(\{x_{j}\})\overset{\mathbb{P}}{{\rightarrow }}+\infty .
\end{equation*}%
Moreover, if $\mathcal{Z}_{D}\neq \varnothing $, then for every $f\in
C_{c}(X)$ satisfying 
\begin{equation*}
f(x_{j})>0\quad \text{for at least one }x_{j}\in \mathcal{Z}_{D},
\end{equation*}%
we have 
\begin{equation*}
\lim_{K\rightarrow +\infty }L_{K}(f)=0.
\end{equation*}

\medskip \noindent \textbf{(iii) Behavior away from the zero set.}

For every $x_{j}\in D_{+}$, we have%
\begin{equation*}
C_{j,K}e^{-K|g(x_{j})|^{2}}\overset{\mathbb{P}}{{\rightarrow }}0.
\end{equation*}%
In particular, 
\begin{equation*}
M_{K}^{(a)}(\{x_{j}\})\overset{\mathbb{P}}{{\rightarrow }}0.
\end{equation*}

\medskip \noindent \textbf{(iv) Asymptotic decomposition.}

The random atomic measure satisfies 
\begin{equation*}
{M_{K}^{(a)}=\sum_{x_{j}\in \mathcal{Z}_{D}}C_{j,K}\delta
_{x_{j}}+\sum_{x_{j}\in D_{+}}C_{j,K}e^{-K|g(x_{j})|^{2}}\delta _{x_{j}}.}
\end{equation*}%
The first component diverges in mass at every true zero, whereas the second
component converges to the zero measure in probability.

Thus, the asymptotic behavior of the driving measure separates the true
zeros from the nonzero candidate locations.
\end{lemma}

\begin{proof}
\textbf{(i) Laplace functional.}

By definition, we can write 
\begin{equation*}
L_{K}(f)=\mathbb{E}\left[ \exp \left\{
-\sum_{j=1}^{m}f(x_{j})C_{j,K}e^{-K|g(x_{j})|^{2}}\right\} \right] .
\end{equation*}%
Using the independence of the Gamma random variables, we obtain 
\begin{equation*}
L_{K}(f)=\prod_{j=1}^{m}\mathbb{E}\left[
e^{-f(x_{j})C_{j,K}e^{-K|g(x_{j})|^{2}}}\right] .
\end{equation*}%
The Laplace transform of a $\func{Gamma}(K,1)$ random variable is 
\begin{equation*}
\mathbb{E}[e^{-tC_{j,K}}]=(1+t)^{-K}.
\end{equation*}%
Substitution gives 
\begin{equation*}
L_{K}(f)=\prod_{j=1}^{m}\left( 1+f(x_{j})e^{-K|g(x_{j})|^{2}}\right) ^{-K},
\end{equation*}%
which establishes the identity.\newline
\medskip \noindent \textbf{(ii) Behavior at true zeros.}\newline
Since 
\begin{equation*}
\mathbb{E}[C_{j,K}]=K,\qquad \func{Var}(C_{j,K})=K,
\end{equation*}%
we have 
\begin{equation*}
\frac{C_{j,K}}{K}\overset{\mathbb{P}}{{\rightarrow }}1,
\end{equation*}%
by Chebyshev's inequality. Hence 
\begin{equation*}
C_{j,K}\overset{\mathbb{P}}{{\rightarrow }}+\infty .
\end{equation*}%
For a true zero, $g(x_{j})=0$, and therefore 
\begin{equation*}
M_{K}^{(a)}(\{x_{j}\})=C_{j,K}.
\end{equation*}%
If $f(x_{j})>0$ at a true zero, the corresponding factor in the Laplace
functional is 
\begin{equation*}
(1+f(x_{j}))^{-K},
\end{equation*}%
which converges to zero. Since all other factors are bounded above by one,
it follows that 
\begin{equation*}
L_{K}(f)\longrightarrow 0.
\end{equation*}%
\noindent \textbf{(iii) Behavior away from the zero set.}\newline
Let 
\begin{equation*}
c_{j}=|g(x_{j})|^{2}>0.
\end{equation*}%
Using Markov's inequality, we have 
\begin{equation*}
\mathbb{P}\left( C_{j,K}e^{-Kc_{j}}>\varepsilon \right) \leq \frac{\mathbb{E}%
[C_{j,K}]e^{-Kc_{j}}}{\varepsilon }=\frac{Ke^{-Kc_{j}}}{\varepsilon }%
\longrightarrow 0.
\end{equation*}%
Thus, 
\begin{equation*}
C_{j,K}e^{-K|g(x_{j})|^{2}}\overset{\mathbb{P}}{{\rightarrow }}0.
\end{equation*}%
\medskip \noindent \textbf{(iv) Asymptotic decomposition.}\newline
The decomposition follows directly from separating the finite index set into
its zero and nonzero subsets. The preceding results establish the respective
limiting behaviors of its two components. This completes the proof.
\end{proof}

\subsection{Statistical consistency}

Next we show that the proposed gamma driven Cox point process is consistent.

\begin{lemma}[Finite-set statistical consistency of the gamma driven Cox
point process]
\label{lem:finite_set_consistency}

Let $g:\mathcal{X}\rightarrow \mathbb{R}$ be a measurable function and let 
\begin{equation*}
D=\{x_{1},\ldots ,x_{m}\}\subset \mathcal{X}
\end{equation*}%
be a finite set of candidate locations.

For each $K>0$, independently generate the atomic random masses 
\begin{equation*}
C_{j,K}\sim \func{Gamma}(K,1),\qquad j=1,\ldots ,m,
\end{equation*}%
and define the random atomic measure 
\begin{equation*}
M_{K}^{(a)}=\sum_{j=1}^{m}C_{j,K}e^{-K|g(x_{j})|^{2}}\delta _{x_{j}}.
\end{equation*}%
Let $M_{K}^{(c)}$ be a diffuse random measure such that 
\begin{equation*}
M_{K}=M_{K}^{(a)}+M_{K}^{(c)}
\end{equation*}%
is locally finite almost surely. Conditional on $M_{K}$, let $N_{K}$ be a
Poisson point process with intensity measure $M_{K}$.

Define the detection indicators 
\begin{equation*}
I_{j,K}=\mathbf{1}\left( N_{K}(\{x_{j}\})\geq 1\right) ,\qquad j=1,\ldots ,m,
\end{equation*}%
and the corresponding true zero indicators 
\begin{equation*}
Z_{j}=\mathbf{1}\left( g(x_{j})=0\right) .
\end{equation*}%
Consider the estimated zero set 
\begin{equation*}
\widehat{Z}_{K}=\{x_{j}\in D:I_{j,K}=1\}.
\end{equation*}%
Then the probability of at least one classification error satisfies 
\begin{equation*}
\mathbb{P}(\widehat{Z}_{K}\neq Z_{D})\leq m_{0}2^{-K}+K\sum_{j:g(x_{j})\neq
0}e^{-K|g(x_{j})|^{2}},
\end{equation*}%
where 
\begin{equation*}
Z_{D}=\{x_{j}\in D:g(x_{j})=0\}
\end{equation*}%
and 
\begin{equation*}
m_{0}=|Z_{D}|,
\end{equation*}%
denotes the number of true zeros in $D$. In particular, 
\begin{equation*}
{\lim_{K\rightarrow +\infty }\mathbb{P}(\widehat{Z}_{K}\neq Z_{D})=0.}
\end{equation*}%
Thus, the proposed gamma driven Cox point process identifies all zeros and
nonzeros correctly on any fixed finite candidate set with probability
approaching one as $K\rightarrow +\infty $.
\end{lemma}

\begin{proof}
A classification error at a candidate location $x_j$ occurs in either of two
ways:

\begin{enumerate}
\item If $g(x_j)=0$, a false negative occurs when 
\begin{equation*}
N_K(\{x_j\})=0.
\end{equation*}

\item If $g(x_j)\neq0$, a false positive occurs when 
\begin{equation*}
N_K(\{x_j\})\geq1.
\end{equation*}
\end{enumerate}

By the gamma driven Cox construction, the detection probability is 
\begin{equation*}
p_{K}(x_{j})=1-\left( 1+e^{-K|g(x_{j})|^{2}}\right) ^{-K}.
\end{equation*}%
Consequently, for a true zero, we have 
\begin{equation*}
\mathbb{P}(I_{j,K}=0)=\left( 1+e^{-K|g(x_{j})|^{2}}\right) ^{-K}=2^{-K},
\end{equation*}%
since $g(x_{j})=0$.

For a nonzero location, the false-positive probability is 
\begin{equation*}
\mathbb{P}(I_{j,K}=1)=1-\left( 1+e^{-K|g(x_{j})|^{2}}\right) ^{-K}\leq
Ke^{-K|g(x_{j})|^{2}},
\end{equation*}%
where we have used the inequality 
\begin{equation*}
1-(1+t)^{-K}\leq Kt,\qquad t\geq 0,\quad K>0.
\end{equation*}%
Let $E_{j,K}$ denote the event of misclassifying $x_{j}$. The event of at
least one classification error is 
\begin{equation*}
\{\widehat{Z}_{K}\neq Z_{D}\}=\bigcup_{j=1}^{m}E_{j,K}.
\end{equation*}%
Applying the union bound gives 
\begin{equation*}
\mathbb{P}(\widehat{Z}_{K}\neq Z_{D})\leq \sum_{j:g(x_{j})=0}\mathbb{P}%
(E_{j,K})+\sum_{j:g(x_{j})\neq 0}\mathbb{P}(E_{j,K})\leq
m_{0}\,2^{-K}+K\sum_{j:g(x_{j})\neq 0}e^{-K|g(x_{j})|^{2}}.
\end{equation*}%
Since $D$ is finite and every nonzero function value satisfies $%
|g(x_{j})|^{2}>0$, each term in the finite sum converges to zero as $%
K\rightarrow +\infty $. Furthermore, 
\begin{equation*}
\lim_{K\rightarrow \infty }m_{0}2^{-K}=0.
\end{equation*}

It follows that 
\begin{equation*}
\lim_{K\rightarrow\infty} \mathbb{P}(\widehat Z_K\neq Z_D)=0.
\end{equation*}

This completes the proof.
\end{proof}

The lemma gives three useful results. First it verifies simultaneous
classification consistency, which is stronger than the pointwise result.
Moreover, it establishes that the entire finite candidate set is classified
correctly with probability approaching one.

Second, it provides an explicit upper bound on classification error. Define 
\begin{equation*}
\Delta :=min_{j:g(x_{j})\neq 0}|g(x_{j})|^{2}.
\end{equation*}%
If there is at least one nonzero candidate, then $\Delta >0$, and the bound
simplifies to 
\begin{equation*}
P(\text{classification error})\leq m_{0}2^{-K}+(m-m_{0})Ke^{-K\Delta }.
\end{equation*}%
This gives a useful way to select $K$ according to a desired error tolerance.

Third, no independence assumption is needed for the error bound. Although
the atomic masses are independent in the construction, the proof uses the
union bound, which does not require independence between the classification
events. The result therefore remains valid under more general dependence
structures, provided the individual detection probabilities retain the
stated bounds.

\subsection{Asymptotics for the binary detection process}

The following theorem establishes three levels of convergence\ for the
binary detection process.

\begin{theorem}[Convergence of the binary zero detection process]
\label{thm:binary_detection_convergence}

Let $\mathcal{X}$ be a locally compact Polish space and let 
\begin{equation*}
D=\{x_{1},\ldots ,x_{m}\}\subset \mathcal{X},
\end{equation*}%
be a finite set of distinct candidate locations. Let $g:\mathcal{X}%
\rightarrow \mathbb{R}$ be measurable, and define the zero set restricted to 
$D$ by 
\begin{equation*}
\mathcal{Z}_{D}=\{x_{j}\in D:g(x_{j})=0\}.
\end{equation*}%
For each $K>0$, let 
\begin{equation*}
C_{j,K}\sim \func{Gamma}(K,1),\qquad j=1,\ldots ,m,
\end{equation*}%
be independent random variables. Define the random atomic measure 
\begin{equation*}
M_{K}=\sum_{j=1}^{m}C_{j,K}e^{-K|g(x_{j})|^{2}}\delta _{x_{j}}.
\end{equation*}%
Conditional on $M_{K}$, let $N_{K}$ be a Poisson point process with
intensity measure $M_{K}$. Define the binary detection process 
\begin{equation*}
\widetilde{N}_{K}:=\sum_{j=1}^{m}\mathbf{1}\{N_{K}(\{x_{j}\})\geq 1\}\delta
_{x_{j}}.
\end{equation*}

Furthermore, define the deterministic zero-counting measure on $D$ by 
\begin{equation*}
N_{\mathcal{Z}}^{D}:=\sum_{x_{j}\in \mathcal{Z}_{D}}\delta _{x_{j}}.
\end{equation*}%
Then, as $K\rightarrow +\infty $, we have%
\begin{equation*}
{\widetilde{N}_{K}}\overset{\mathbb{P}}{{\rightarrow }}{N_{Z}^{D}.}
\end{equation*}%
More precisely, for every continuous function $f:X\rightarrow \mathbb{R}$,
we can write 
\begin{equation*}
\int\limits_{\mathcal{X}}{f(x)\widetilde{N}_{K}(dx)\overset{\mathbb{P}}{{%
\rightarrow }}\int\limits_{\mathcal{X}}f(x)N_{\mathcal{Z}}^{D}(dx),}
\end{equation*}%
and in particular, we have%
\begin{equation*}
\widetilde{N}_{K}(D)\overset{\mathbb{P}}{{\rightarrow }}|\mathcal{Z}_{D}|.
\end{equation*}%
Thus, the binary detection process converges in probability to the counting
measure supported on the true zeros among the candidate locations.
\end{theorem}

\begin{proof}
For each candidate location $x_{j}$, define 
\begin{equation*}
I_{j,K}=\mathbf{1}\{N_{K}(\{x_{j}\})\geq 1\}
\end{equation*}%
and 
\begin{equation*}
Z_{j}=\mathbf{1}\{g(x_{j})=0\}.
\end{equation*}%
By the gamma driven Cox construction, 
\begin{equation*}
\mathbb{P}(I_{j,K}=1)=1-\left( 1+e^{-K|g(x_{j})|^{2}}\right) ^{-K}.
\end{equation*}%
The asymptotic detection property gives 
\begin{equation*}
\mathbb{P}(I_{j,K}=1)\longrightarrow Z_{j}.
\end{equation*}%
Since $I_{j,K}$ and $Z_{j}$ are binary random variables, we obtain 
\begin{equation*}
\mathbb{P}(I_{j,K}\neq Z_{j})=%
\begin{cases}
\mathbb{P}(I_{j,K}=0), & Z_{j}=1, \\ 
\mathbb{P}(I_{j,K}=1), & Z_{j}=0,%
\end{cases}%
\longrightarrow 0.
\end{equation*}%
Now, because $D$ is finite, the union bound yields 
\begin{equation*}
\mathbb{P}\left( \widetilde{N}_{K}\neq N_{\mathcal{Z}}^{D}\right) =\mathbb{P}%
\left( \bigcup_{j=1}^{m}\{I_{j,K}\neq Z_{j}\}\right) \leq \sum_{j=1}^{m}%
\mathbb{P}(I_{j,K}\neq Z_{j})\longrightarrow 0.
\end{equation*}%
Therefore, 
\begin{equation*}
\widetilde{N}_{K}\overset{\mathbb{P}}{{\rightarrow }}N_{\mathcal{Z}}^{D}.
\end{equation*}%
To establish the integral convergence, observe that on the event $\{%
\widetilde{N}_{K}=N_{\mathcal{Z}}^{D}\}$, we have 
\begin{equation*}
\int\limits_{\mathcal{X}}f(x)\widetilde{N}_{K}(dx)=\int\limits_{\mathcal{X}%
}f(x)N_{\mathcal{Z}}^{D}(dx).
\end{equation*}%
Since the probability of the complementary event converges to zero, the
integrals converge in probability.

Finally, choosing $f\equiv 1$ on the finite set $D$ gives 
\begin{equation*}
\widetilde{N}_{K}(D)\overset{\mathbb{P}}{{\rightarrow }}|\mathcal{Z}_{D}|,
\end{equation*}%
which completes the proof.
\end{proof}

There are three main points to consider from the latter theorem. First, the
individual locations; each candidate is classified correctly with
probability approaching one. Second the probability of recovering every zero
in the finite candidate set approaches one. Third the entire binary
detection process converges in probability to the deterministic zero
counting measure.

This third result is particularly interesting from the perspective of point
process theory. Along with the lemma on the asymptotic behavior of the
driving random measures, these results establish a logical progression,
where the driving measure concentrates its mass at true zeros, the Cox
process generates increasingly many points at those locations, and the
binary detection process converges to the zero counting measure. Therefore,
the proposed zero detection construction is integrally linked with the
mathematical theory of random measures and stochastic point processes.

\section{Zero detection for Riemann's $\protect\xi $\ function}

Since the theoretical development is agnostic to the event space $\mathcal{X}
$, we will discuss an important example on the complex plane $\mathbb{C}$.

Consider Riemann's $\zeta $ function (\cite{riemann1859ueber}, see \cite%
{selberg1946contribution}, \cite{edwards1974riemannaes}, \cite%
{odlyzko1987distribution}, \cite{patterson1995introduction}, \cite%
{borwein2000efficient}, \cite{conrey2003riemann}, \cite{bombieri2006riemann}%
, \cite{guillera2008double}, \cite{hiary2012zeta}, \cite{bober2018new}, \cite%
{pratt2020more}, \cite{spigler2025brief}, \cite{bellotti2025explicit})
defined for a complex number $s=\sigma +it\in \mathbb{C},$ with $\sigma >1,$
by%
\begin{equation}
\zeta (s):=\sum\limits_{n=1}^{+\infty }\frac{1}{n^{s}}.  \label{Zetafunc1}
\end{equation}%
The $\zeta $-function is of paramount importance in number theory, since it
provides insight on the distribution of prime numbers.

Recall that the $\zeta $-function can be defined to be meromorphic in the
plane, i.e., an analytic function except at its poles, with only a simple
pole at $s=1$. For a rigorous exposition of this construction see \cite%
{conway2012functions} (Theorem 8.14, pg 189). More precisely, for $s\neq 1,$ 
$\zeta $ satisfies Riemann's functional equation%
\begin{equation}
\zeta (s)=2(2\pi )^{s-1}\Gamma (1-s)\zeta (1-s)\sin \left( \frac{\pi s}{2}%
\right) ,  \label{functionaleq}
\end{equation}%
where $\Gamma $ denotes the gamma function, define for $\func{Re}(s)>0$ by%
\begin{equation}
\Gamma (s):=\int\limits_{0}^{+\infty }x^{s-1}e^{-x}dx.  \label{gammafunc}
\end{equation}%
The meromorphic continuation of $\zeta $ is holomorphic in $\func{Re}(s)<0$;
its zeros in that region are precisely the trivial zeros at the negative
even integers. The functional equation, together with the pole structure of $%
\Gamma $ and the nonvanishing of $\zeta (1-s)$ for $\func{Re}(1-s)>1$, shows
that the zeros in $\func{Re}(s)<0$ are $s=-2,-4,-6,...$. These are called
the trivial zeros and any other zeros are called nontrivial zeros.

Riemann's hypothesis (RH, see \cite{riemann1859ueber}) is a conjecture
asserting that every nontrivial zero has real part $0.5$, i.e., the
non-trivial zeros are obtained on the critical line $0.5+it$, for some $t\in 
\mathbb{R}$. The first few positive imaginary parts (with their negatives
also zeros by conjugation) that have been discovered include $0.5\pm
14.134725...i$, $0.5\pm 21.022040...i,$ $0.5\pm 25.010858...i,$ $0.5\pm
30.424876...i,$ $0.5\pm 32.935062...i,$ $0.5\pm 37.586178...i,$ and $0.5\pm
40.918719...i.$

Furthermore, except for the trivial zeros, $\zeta $ has no other zeros
outside the open strip%
\begin{equation}
\mathcal{S}:=\{s=\sigma +it\in \mathbb{C}:0<\sigma <1,t\in \mathbb{R}\},
\label{strip1}
\end{equation}%
known as the critical strip. We use the term open critical strip for $0<%
\func{Re}(s)<1$. The zeta function has no zeros on either boundary line $%
\func{Re}(s)=0$ or $\func{Re}(s)=1$, and thus every nontrivial zero lies in
this open strip. Here we dropped the lines for $\sigma =0$, and $1$, from
the standard definition of $\mathcal{S}$ (see \cite{conway2012functions},
Definition 8.16, pg 189)$,$ in order to have $\mathcal{S}$ open, which is
not crucial in studying the zeros of $\zeta $, since it has been shown (see 
\cite{hadamard1896distribution} and \cite{poussin1897recherches}), that
there are no zeros of $\zeta $ on the line $\sigma =1,$ and hence by (\ref%
{functionaleq}), no zeros on $\sigma =0$. In addition, it has been
established (see \cite{hardy1914zeros}) that there are infinitely many zeros
on the critical line $\sigma =\frac{1}{2}$. Since the complex conjugate of $%
\zeta $ satisfies $\zeta (\overline{s})=\overline{\zeta (s)},$ then if $%
\kappa $ is a non-trivial zero of $\zeta ,$ we must have that $\overline{%
\kappa },$ $1-\kappa ,$ and $1-\overline{\kappa }$ are also non-trivial
zeros of $\zeta $, so that instead of $t\in \mathbb{R}$ in (\ref{strip1}), a
lot of attention has been placed on the open half-strip, for $t>0.$ The
functional equation gives the reflection $\kappa \mapsto 1-\kappa ,$ and
conjugation gives the other symmetries. A last consideration involves the
multiplicity of the zeros$;$ all the zeros of $\zeta (s)$ computed up to a
specific height in (see \cite{van1986zeros} and the discussion in \cite%
{bombieri2006riemann} and the references therein), are simple.

The completed zeta function is Riemann's $\xi $ function (see \cite%
{riemann1859ueber}) which is an entire function and is defined by%
\begin{equation}
\xi (s):=\frac{1}{2}s(s-1)\pi ^{-\frac{s}{2}}\Gamma \left( \frac{s}{2}%
\right) \zeta (s)=(s-1)\pi ^{-\frac{s}{2}}\Gamma \left( \frac{s}{2}+1\right)
\zeta (s),  \label{ksifunc}
\end{equation}%
for all $s\in \mathbb{C},$ with the second equality in (\ref{ksifunc}), the
now standard notation for $\xi ,$ which follows from Landau (see \cite%
{edwards1974riemannaes}, pg 16). Then $\xi (s)$ is an entire function (i.e.,
a function which is defined and analytic in the whole complex plane $\mathbb{%
C}$) of order one, it satisfies the symmetries 
\begin{equation}
\xi (s)=\xi (1-s),\text{ and }\xi (\overline{s})=\overline{\xi (s)},
\label{symmetries}
\end{equation}%
and its zeros, counted with multiplicity, are precisely the non-trivial
zeros of $\zeta (s)$. In particular, all of its zeros lie in $\mathcal{S}$;
specifically, $s=0,1$ are not zeros of $\xi $ as one might be tempted to
deduce from the first equality of equation (\ref{ksifunc})\ (in fact $\xi
(0)=\xi (1)=\frac{1}{2}),$ since the pole of the $\zeta $-function at $s=1$,
is canceled out by the term $(s-1)$. Moreover, $\Gamma \left( s\right) $ is
meromorphic with simple poles at the integers $s=0,-1,-2,...$ (see \cite%
{conway2012functions}, Definition 7.2, pg 172)$,$ such that, these poles are
canceled by the trivial zeros of $\zeta .$

\subsection{The real-valued reduction}

The basic Cox theorem is formulated for a real-valued measurable function.
There are two natural reductions.

\begin{definition}[Complex-plane detection function]
Define 
\begin{equation}
g_{\mathbb{C}}(s):=|\xi (s)|,\qquad s\in \mathbb{C}.  \label{eq:gcomplex}
\end{equation}%
Then $g_{\mathbb{C}}:\rightarrow \lbrack 0,\infty )$ is continuous and 
\begin{equation}
g_{\mathbb{C}}(s)=0\Leftrightarrow \xi (s)=0.  \label{eq:zeroequiv}
\end{equation}
\end{definition}

Continuity follows because $\xi $ is entire and the complex modulus is
continuous. The function $g_{\mathbb{C}}$ is generally not holomorphic; this
is irrelevant to the Cox construction, which requires measurability rather
than analyticity of the detection function.

There is also a real-valued restriction to the critical line. Define 
\begin{equation}
h(t):=\xi \left( \frac{1}{2}+it\right) ,\qquad t\in \mathbb{R}.
\label{eq:hcritical}
\end{equation}%
By (\ref{symmetries}), we have 
\begin{equation*}
\overline{h(t)}=\overline{\xi (1/2+it)}=\xi (1/2-it)=\xi (1-(1/2+it))=\xi
(1/2+it)=h(t).
\end{equation*}%
Thus $h(t)\in \mathbb{R}$ for every real $t$, and 
\begin{equation}
h(t)=0\Leftrightarrow \xi (1/2+it)=0.  \label{eq:hzeros}
\end{equation}

For completeness, let 
\begin{equation*}
\vartheta (t)=\func{Im}\log \Gamma \left( \frac{1}{4}+\frac{it}{2}\right) -%
\frac{t}{2}\log \pi ,
\end{equation*}%
with a continuous choice of argument along the line, and define Hardy's
function 
\begin{equation}
Z(t)=e^{i\vartheta (t)}\zeta \left( \frac{1}{2}+it\right) .
\label{eq:hardyz}
\end{equation}%
The standard conjugation symmetry implies that $Z(t)$ is real for real $t$.
Direct substitution into (\ref{ksifunc}) gives 
\begin{equation}
\xi \left( \frac{1}{2}+it\right) =-\frac{1}{2}\left( t^{2}+\frac{1}{4}%
\right) \pi ^{-1/4}\left\vert \Gamma \left( \frac{1}{4}+\frac{it}{2}\right)
\right\vert Z(t).  \label{eq:xizrelation}
\end{equation}%
The factor multiplying $Z(t)$ is real and nonzero for every real $t$.
Therefore 
\begin{equation}
Z(t)=0\Leftrightarrow h(t)=0\Leftrightarrow \xi (1/2+it)=0.
\label{eq:Zzeros}
\end{equation}%
Either $h$ or $Z$ can be used as the real-valued function in a
one-dimensional Cox detector.

\subsection{Critical line detector}

We now specialize to candidates on the critical line. Let 
\begin{equation*}
\mathcal{T}=\{t_{1},\ldots ,t_{m}\}\subset \mathbb{R},
\end{equation*}%
be a finite set of distinct real ordinates, and let 
\begin{equation*}
D_{\mathrm{crit}}=\left\{ \frac{1}{2}+it_{j}:1\leq j\leq m\right\} .
\end{equation*}%
We will use the real-valued function $h(t)=\xi (1/2+it)$ from (\ref%
{eq:hcritical}). Alternatively, use Hardy's real-valued $Z(t)$. For
definiteness, write $q(t)$ for either $h(t)$ or $Z(t)$; both have the same
real zero set by (\ref{eq:Zzeros}).

For independent $C_{j,K}\sim \func{Gamma}(K,1)$, define 
\begin{equation}
M_{K,\mathrm{crit}}^{(a)}=\sum_{j=1}^{m}C_{j,K}e^{-K|q(t_{j})|^{2}}\delta
_{t_{j}}.  \label{eq:critical-measure}
\end{equation}%
Conditional on this measure (optionally augmented by a diffuse component),
let $N_{K,\mathrm{crit}}$ be a Poisson point process with this intensity.
Define 
\begin{equation*}
I_{j,K}^{\mathrm{crit}}=\mathbf{1}\{N_{K,\mathrm{crit}}(\{t_{j}\})\geq 1\}.
\end{equation*}%
Exactly the same Laplace transform calculation as in Lemma \ref%
{lem:random_measure_asymptotics}\ yields 
\begin{equation}
\mathbb{P}(I_{j,K}^{\mathrm{crit}}=1)=1-\left( 1+e^{-K|q(t_{j})|^{2}}\right)
^{-K}.  \label{eq:critical-prob}
\end{equation}%
Consequently, we can write 
\begin{equation}
\lim_{K\rightarrow \infty }\mathbb{P}(I_{j,K}^{\mathrm{crit}}=1)=\mathbf{1}%
\{q(t_{j})=0\}=\mathbf{1}\{\xi (1/2+it_{j})=0\}.  \label{eq:critical-limit}
\end{equation}%
For the finite set $\mathcal{T}$, if 
\begin{equation*}
T_{0}=\{t_{j}:q(t_{j})=0\},\qquad m_{0}=|T_{0}|,
\end{equation*}%
then 
\begin{equation}
\mathbb{P}(\widehat{T}_{K}\neq T_{0})\leq
m_{0}\,2^{-K}+K\sum_{j:q(t_{j})\neq 0}e^{-K|q(t_{j})|^{2}},
\label{eq:critical-error}
\end{equation}%
where $\widehat{T}_{K}=\{t_{j}:I_{j,K}^{\mathrm{crit}}=1\}$. Thus the method
consistently classifies all ordinates in any fixed finite list.

A natural question that arises is which real-valued function should be used?
The choice $q=h$ gives the direct restriction of the completed zeta function
to the critical line. The choice $q=Z$ removes the nonzero gamma and
elementary factors in (\ref{eq:xizrelation}); it has the same zeros and is a
standard real-valued representation for computations on the critical line.
The numerical scale of $h$ and $Z$ differs, so the effective finite $K$
behavior differs unless $K$ is adjusted accordingly. Their asymptotic zero
classification statements are identical.

\subsection{A three point example}

For numerical convenience, introduce a fixed rescaling 
\begin{equation*}
g(t)=10^{4}h(t).
\end{equation*}%
Multiplication by a nonzero constant leaves the zero set unchanged. It does,
however, change the numerical values in the exponential weight, so the
choice of scale and the parameter $K$ must be considered together.

Let 
\begin{equation*}
D=\{t_{1},t_{2},t_{3}\},
\end{equation*}%
where 
\begin{align*}
t_1&=14.1347251417346937904572519836\ldots, \\
t_2&=14, \\
t_3&=14.5.
\end{align*}%
Here $t_{1}$ is the ordinate of the first nontrivial zero, used as a known
positive control. The other two points are negative controls. High precision
evaluation gives approximately 
\begin{equation*}
\begin{array}{c|c|c}
\toprule\text{Candidate} & |h(t_{j})| & |g(t_{j})|^{2} \\ 
\midrule t_{1} & \text{numerically zero} & \text{numerically zero} \\ 
14 & 2.01294444235\times 10^{-4} & 4.0519453280 \\ 
14.5 & 4.06827079891\times 10^{-4} & 16.5508272933 \\ 
\bottomrule &  & 
\end{array}%
\end{equation*}%
The first row should be understood with care; when the zero ordinate is
computed at high precision, the residual is governed by numerical precision.
It is not mathematically correct to interpret a floating-point residual as
the exact value of the function at the exact zero.

\subsection{The gamma driven Cox construction}

Fix 
\begin{equation*}
K=10.
\end{equation*}%
For each $t_{j}\in D$, let 
\begin{equation*}
C_{j,10}\overset{\mathrm{ind}}{\sim }\func{Gamma}(10,1),
\end{equation*}%
where the parameters denote shape and scale. Thus 
\begin{equation*}
f_{C}(c)=\frac{c^{9}e^{-c}}{\Gamma (10)},\qquad c>0,\qquad \mathbb{E}\left(
C_{j,10}\right) =10.
\end{equation*}%
Define the random atomic measure 
\begin{equation*}
M_{10}=\sum_{j=1}^{3}C_{j,10}\exp \left( -10|g(t_{j})|^{2}\right) \delta
_{t_{j}}.
\end{equation*}%
Conditional on $M_{10}$, let $N$ be a Poisson point process with intensity
measure $M_{10}$. Since $M_{10}$ is supported on three atoms, this is
equivalent to generating conditionally independent random variables 
\begin{equation*}
N_{j}\mid C_{j,10}\sim \func{Poisson}\left(
C_{j,10}e^{-10|g(t_{j})|^{2}}\right) ,\qquad j=1,2,3.
\end{equation*}%
The detector declares $t_{j}$ detected precisely when 
\begin{equation*}
\widehat{Z}_{j}=\mathbf{1}_{\{N_{j}\geq 1\}}=1.
\end{equation*}%
The Poisson counts must be generated conditionally on the gamma variables.
Replacing this two stage construction by Poisson variables with only the
unconditional mean would not implement the same Cox process.

\subsection{Exact detection probability and implementation}

Conditional on $C_{j,K}$, 
\begin{equation*}
\mathbb{P}(N_{j}=0\mid C_{j,K})=\exp \left(
-C_{j,K}e^{-K|g(t_{j})|^{2}}\right) .
\end{equation*}%
For $C\sim \func{Gamma}(K,1)$, its Laplace transform is 
\begin{equation*}
\mathbb{E}(e^{-uC})=(1+u)^{-K},\qquad u\geq 0.
\end{equation*}%
Taking expectations therefore gives 
\begin{equation*}
\mathbb{P}(N_{j}=0)=\left( 1+e^{-K|g(t_{j})|^{2}}\right) ^{-K}.
\end{equation*}%
Hence 
\begin{equation*}
p_{j,K}:=\mathbb{P}(N_{j}\geq 1)=1-\left( 1+e^{-K|g(t_{j})|^{2}}\right)
^{-K}.
\end{equation*}

For an exact zero, $g(t_j)=0$, and 
\begin{equation*}
p_{j,K}=1-2^{-K}.
\end{equation*}
For a nonzero candidate, the elementary inequality $1-(1+x)^{-K}\leq Kx$,
for $x\geq0$, yields 
\begin{equation*}
p_{j,K}\leq K e^{-K|g(t_j)|^2}.
\end{equation*}

With $K=10$, the known zero has detection probability 
\begin{equation*}
p_{1,10}=1-2^{-10}=0.9990234375,
\end{equation*}%
in the idealized exact-zero calculation. For $t_{2}=14$, 
\begin{equation*}
e^{-10|g(14)|^{2}}=e^{-40.51945328\ldots },
\end{equation*}%
and consequently 
\begin{equation*}
p_{2,10}\approx 2.5271\times 10^{-17}.
\end{equation*}%
For $t_{3}=14.5$, 
\begin{equation*}
e^{-10|g(14.5)|^{2}}=e^{-165.50827293\ldots },
\end{equation*}%
so that 
\begin{equation*}
p_{3,10}\approx 1.3203\times 10^{-71}.
\end{equation*}%
Thus the idealized single run probabilities are 
\begin{equation*}
\begin{array}{c|c|c}
\toprule\text{Candidate} & \text{Detection probability} & \text{Expected
classification} \\ 
\midrule14.134725\ldots  & 0.9990234375 & \text{detected} \\ 
14 & 2.5271\times 10^{-17} & \text{not detected} \\ 
14.5 & 1.3203\times 10^{-71} & \text{not detected} \\ 
\bottomrule &  & 
\end{array}%
\end{equation*}%
These are theoretical probabilities, not guaranteed outcomes of a particular
random realization. A true zero can be missed with probability $2^{-K}$,
while a nonzero point can be falsely detected with a small positive
probability.

\subsection{Finite-set consistency theorem}

Let $D=\{t_1,\ldots,t_m\}$ be a fixed finite set, and let 
\begin{equation*}
Z_D=\{t_j\in D:h(t_j)=0\}
\end{equation*}
be its true zero subset. Let $\widehat Z_K$ be the subset detected by the
Cox construction with parameter $K$.

\begin{proposition}
For every $K>0$, 
\begin{equation*}
\mathbb{P}(\widehat{Z}_{K}\neq Z_{D})\leq |Z_{D}|2^{-K}+K\sum_{t_{j}\in
D\setminus Z_{D}}e^{-K|g(t_{j})|^{2}}.
\end{equation*}%
In particular, for every fixed finite candidate set $D$, we have 
\begin{equation*}
\mathbb{P}(\widehat{Z}_{K}\neq Z_{D})\longrightarrow 0\qquad \text{as }%
K\longrightarrow +\infty .
\end{equation*}
\end{proposition}

\begin{proof}
An error occurs if at least one true zero is missed or at least one nonzero
candidate is detected. By the union bound, we have 
\begin{equation*}
\mathbb{P}(\widehat{Z}_{K}\neq Z_{D})\leq \sum_{t_{j}\in
Z_{D}}(N_{j}=0)+\sum_{t_{j}\in D\setminus Z_{D}}(N_{j}\geq 1).
\end{equation*}%
For $t_{j}\in Z_{D}$, $g(t_{j})=0$, so 
\begin{equation*}
\mathbb{P}(N_{j}=0)=2^{-K}.
\end{equation*}%
For $t_{j}\notin Z_{D}$, the detection probability formula and $%
1-(1+x)^{-K}\leq Kx$ give 
\begin{equation*}
\mathbb{P}(N_{j}\geq 1)\leq Ke^{-K|g(t_{j})|^{2}}.
\end{equation*}%
Combining these inequalities proves the bound. Since $D$ is finite and $%
g(t_{j})\neq 0$ for every $t_{j}\in D\setminus Z_{D}$, each term in the
second sum tends to zero as $K\rightarrow +\infty $. The first sum also
tends to zero. This proves the asserted convergence.
\end{proof}

Equivalently, if $D\setminus Z_{D}\neq \varnothing $, define 
\begin{equation*}
\Delta _{D}=\min_{t_{j}\in D\setminus Z_{D}}|g(t_{j})|^{2}>0.
\end{equation*}%
Then 
\begin{equation*}
\mathbb{P}(\widehat{Z}_{K}\neq Z_{D})\leq \left\vert Z_{D}\right\vert
2^{-K}+(m-|Z_{D}|)Ke^{-K\Delta _{D}}\longrightarrow 0.
\end{equation*}

\subsection{Discussion}

The method does not distinguish an exact zero from an arbitrarily close
nonzero point with certainty at finite $K$. If $|g(t)|$ is small, then $%
e^{-K|g(t)|^{2}}$ is close to one, and the detection probability is close to 
$1-2^{-K}$, whether or not $g(t)$ is exactly zero. The rescaling $g=10^{4}h$
is not innocuous for a fixed numerical value of $K$; replacing $h$ by $ag$
changes the weight from $e^{-K|h|^{2}}$ to $e^{-K|ah|^{2}}$. The zero set is
invariant, but the finite $K$ operating characteristics are not.

The binary detector records whether a candidate atom receives at least one
point. It does not determine the multiplicity of a zero, and it does not
prove that the candidate list is complete. For rigorous computational
claims, one can use the Cox detector as a screening device and combine it
with validated zero isolation or the argument principle. This approach will
be pursued elsewhere.

The example gives an explicit implementation of the gamma driven Cox process
method for $h(t)=\xi (1/2+it)$. At each candidate ordinate, the algorithm
draws a gamma random intensity multiplier and then a conditional Poisson
count, and the exact detection probability follows from the gamma Laplace
transform. On a fixed finite candidate set, the probability of classifying
every candidate correctly tends to one as $K$ increases.

Note that this is a rigorous finite set consistency statement, not a proof
of the Riemann Hypothesis and not a guarantee that a numerical grid contains
every zero. The method is best viewed as a probabilistic screening procedure
whose output should be followed by independent mathematical certification.

\section{Concluding remarks}

The construction of the gamma driven Cox point process $N_{K}$ combines two
sources of information. The deterministic factor $\exp \{-K|g(x_{j})|^{2}\}$
assigns greater weight to locations where the absolute function (or modulus)
value is small, while the gamma distributed random mass introduces
stochastic variability into the corresponding atomic intensity. The
parameter $K$ controls the asymptotic behavior of the construction.

In particular, the probability of detecting at least one point at a
candidate location $x_{j}$ can be derived explicitly. Writing 
\begin{equation*}
I_{j,K}=1\{N_{K}(\{x_{j}\})\geq 1\},
\end{equation*}%
the gamma Laplace transform gives 
\begin{equation*}
\mathbb{P}(I_{j,K}=1)=1-\mathbb{E}\left[ \exp \left\{ -C_{j,K}\exp
(-K|g(x_{j})|^{2})\right\} \right] =1-\left( 1+\exp
(-K|g(x_{j})|^{2})\right) ^{-K}.
\end{equation*}%
This expression establishes the fundamental detection property of the
proposed model; namely, for every fixed candidate location $x_{j}$, we have 
\begin{equation*}
\lim_{K\rightarrow \infty }\mathbb{P}(I_{j,K}=1)=1\{g(x_{j})=0\}.
\end{equation*}%
Thus, the detection probability converges to one at a true zero and to zero
at a nonzero location.

The original driving measure has diverging mass at true zeros, and thus, the
limiting object is not an ordinary locally finite random measure. The binary
detection process solves this problem for finite candidate sets, but an
extension to continuous domains will require a different construction or a
suitable notion of convergence.

Note that the prescribed finite set of candidate locations needs to be
supplied in the methods presented herein, and then in our theoretical
results we have established a rigorous criterion for zero detection. The
other direction is of particular interest; can we use the created PPPP or
CPP to find the actual value of a zero of $g$ using statistical simulation?
The answer is a definite yes, and the theory, along with applications of
this idea are currently being developed. These results are of great interest
and will be presented elsewhere.

\section*{Declarations}

\begin{itemize}
\item I have no conflict of interest with the results in this paper.

\item This research received no funding.
\end{itemize}

\bibliographystyle{plainnat}
\bibliography{ZetaPP}

\end{document}